\documentclass[10pt]{amsart}
\usepackage{amssymb,amsmath,amscd,amsfonts}
\input diagxy

\newcommand{\Ker}{\mathop{\rm Ker}}
\newcommand{\Img}{\mathop{\rm Im}}
\newcommand{\Coker}{\mathop{\rm Coker}}
\newcommand{\coker}{\mathop{\rm coker}}
\newcommand{\im}{\mathop{\rm im}}
\newcommand{\coim}{\mathop{\rm coim}}
\newcommand{\Coim}{\mathop{\rm Coim}}
\newcommand{\id}{\mathop{\rm id}}
\newtheorem{axiom}{Axiom}
\newtheorem{them}{Theorem}
\newtheorem{corl}{Corollary}
\newtheorem{lemma}{Lemma}
\newtheorem{remk}{Remark}

\begin{document}
\title[Lambek Invariants in a Quasi-Abelian Category]
{On the Lambek Invariants of Commutative Squares in
a~Quasi-Abelian Category}
\author{Yaroslav Kopylov}
\address{Yaroslav Kopylov
\newline\hphantom{iii} Sobolev Institute of Mathematics,
\newline\hphantom{iii} Akademik Koptyug Prospect 4,
\newline\hphantom{iii} 630090, Novosibirsk, Russia}%
\email{yakop@math.nsc.ru}
%\classification{Mathematics Subject Classifications (2000)}{18A20, 46M99}
\maketitle
\begin{abstract}
We consider the invariants $\Ker$ and $\Img$ for commutative squares
in quasi-abelian categories. These invariants were introduced
by Lambek for groups and then studied by Hilton and Nomura in exact
categories.

\vspace{2mm}
\noindent
\textbf{Key words and phrases:} quasi-abelian category, commutative square, Lambek invariants
\end{abstract}

\section{Introduction}
\label{intro}
In 1964, Lambek introduced the following invariants for a~commutative
square
\begin{equation}\label{1}
\bfig
\morphism(0,200)<400,0>[C`D;\alpha]
\morphism(0,200)<0,-400>[C`A;g]
\morphism(0,-200)<400,0>[A`B;\beta]
\morphism(400,200)<0,-400>[D`B;f]
\place(200,0)[S]
\efig
\end{equation}
in the category of groups:
$$
\Img S=(\Img\beta\cap\Img f)/\Img(f\alpha), \quad
\Ker S=\Ker(f\alpha)/(\Ker\alpha+\Ker g).
$$
In~\protect\cite{La}, he proved the following assertion.
\smallskip

\textit{
Given a~commutative diagram
\begin{equation}\label{2}
\bfig
\morphism(0,200)<400,0>[A`B;f]
\morphism(400,200)<400,0>[B`C;g]
\morphism(0,200)<0,-400>[A`A';a]
\morphism(400,200)<0,-400>[B`B';b]
\morphism(800,200)<0,-400>[C`C';c]
\morphism(0,-200)<400,0>[A'`B';f']
\morphism(400,-200)<400,0>[B'`C';g']
\place(200,0)[S]
\place(600,0)[T]
\efig
\end{equation}
of groups and group homomorphisms with exact rows, there is a~natural
isomorphism
$$
\Lambda:\Img S\stackrel{\cong}{\longrightarrow} \Ker T.
$$ }
\smallskip

Later Leicht extended this theorem to arbitrary exact categories
(see~\protect\cite{Le}). In~\protect\cite{No1,No2},
Nomura considered the case where the rows in~(\ref{2}) are not exact
but only semiexact, constructed a canonical morphism
$\Lambda:\Img S\to\Ker T$, and proved that there is an exact sequence
\begin{multline}\label{3}
0\rightarrow H(\Ker(bf)\rightarrow\Ker b\rightarrow \Ker c)
\rightarrow \Ker(H\rightarrow H') \rightarrow
\Img S \stackrel\Lambda\rightarrow \Ker T \\
\rightarrow \Coker(H\rightarrow H')
\rightarrow H(\Coker a\rightarrow \Coker b \rightarrow \Coker(g'b))
\rightarrow 0,
\end{multline}
where the arrows between the kernels and cokernels in parentheses are
natural morphisms, $H(\cdot\rightarrow\cdot\rightarrow\cdot)$ stands for
the homology of the $0$-sequence in parentheses,
$H=H(A\rightarrow B\rightarrow C)$, and $H'=H(A'\rightarrow B'\rightarrow C')$.

In this paper, we study the Lambek invariants in quasi-abelian
categories, first considered by Ra{\u\i}kov in~\protect\cite{Ra} under
the name of semiabelian categories. Apart from all abelian categories, the class
of quasi-abelian categories contains many nonabelian additive
categories of functional analysis and topological algebra. The categories of (Hausdorff or all) topological
abelian groups, topological vector spaces, Banach (or normed) spaces, filtered modules over
filtered rings, and torsion-free abelian groups are typical examples of quasi-abelian
categories. The main difference between the quasi-abelian and abelian categories lies in
the fact that the standard diagram lemmas hold in quasi-abelian categories under some
extra conditions which usually amount to the strictness of some morphisms.
Quasi-abelian categories have been actively studied in the recent
years (see~\protect\cite{GK1,GK2,K,KK1,KK2,P1,P2,P3,R1,R2,R3,S}).

In the category $\mathcal{B}an$ of Banach spaces
topological abelian groups, the strictness of a morphism $\alpha$ means that
the range of $\alpha$ is closed. In the category of topological abelian
groups, a morphism $\alpha$ strict if and only if its image is closed
and, moreover, $\alpha$ maps open sets onto open sets.

In a~quasi-abelian category, Nomura's morphism
$\Lambda:\Img S\to\Ker T$ is
defined only if $b$ is strict in~(\ref{2}) because the definition uses
the fact that $b$ is the composition of its image and coimage. Lambek's
isomorphism holds under the same condition (see~\protect\cite{No1}).

The structure of the paper is as follows. In Section~{\ref{sec1}},
we recall some basic definitions and facts about quasi-abelian
categories. In Section~{\ref{sec2}},
we construct a~morphism $\zeta:\Ker T\to\Img S$ for a~diagram~(\ref{2}) with
exact rows in the general case and suggest quasi-abelian versions for some
assertions proved by Nomura~\protect\cite{No1} and Hilton~\protect\cite{Hi}
for abelian and exact categories.

\section{Quasi-Abelian Categories}\label{sec1}

We consider additive categories satisfying the following axiom.

\begin{axiom}\label{a1}
Each morphism has kernel and cokernel.
\end{axiom}

We denote by $\ker\alpha$ ($\coker\alpha$) an arbitrary kernel (cokernel)
of~$\alpha$ and by $\Ker\alpha$ ($\Coker\alpha$)
the corresponding object; the equality $a=\ker b$ ($a=\coker b$) means
that $a$ is a kernel of~$b$ ($a$ is a~cokernel of~$b$).

In a category meeting Axiom~\ref{a1}, every morphism
$\alpha$
admits a canonical decomposition
$\alpha=(\im\alpha)\overline\alpha(\coim\alpha)=(\im\alpha)\widetilde\alpha$,
where
$\im\alpha=\ker\coker\alpha$,
$\coim\alpha=\coker\ker\alpha$. Two canonical decompositions of the same
morphism are obviously naturally isomorphic.
A morphism
$\alpha$
is called {\it strict} if
$\overline\alpha$
is an isomorphism.

We use the following notations of~\protect\cite{KC}:

$O_c$ is the class of all strict morphisms,

$M$ is the class of all monomorphisms,

$M_c$ is the class of all strict monomorphisms ($=$ kernels),

$P$ is the class of all epimorphisms,

$P_c$ is the class of all strict epimorphisms ($=$ cokernels).

\begin{lemma}[\protect\cite{BD,EH1,KC,Ra}]\label{l1}
The following assertions hold in an additive category meeting
Axiom~{\rm\ref{a1}:}

{\rm(1)} $\ker\alpha\in M_c$ and $\coker\alpha\in P_c$
for every
$\alpha${\rm;}

{\rm(2)} $\alpha\in M_c \Longleftrightarrow \alpha=\im\alpha$,
$\alpha\in P_c \Longleftrightarrow \alpha=\coim\alpha${\rm;}

{\rm(3)} a morphism
$\alpha$
is strict if and only if it is representable in the form
$\alpha=\alpha_1 \alpha_0$
with
$\alpha_0\in P_c$,
$\alpha_1\in M_c${\rm;}
in every such representation,
$\alpha_0=\coim\alpha$
and
$\alpha_1=\im\alpha${\rm;}

{\rm(4)} if some commutative square
\begin{equation*}
\bfig
\morphism(0,200)<400,0>[C`D;\alpha]
\morphism(0,200)<0,-400>[C`A;g]
\morphism(0,-200)<400,0>[A`B;\beta]
\morphism(400,200)<0,-400>[D`B;f]
\efig
\end{equation*}
is a pullback then $\ker f=\alpha(\ker g)$ and $f=\ker\xi$ implies $g=\ker(\xi\beta)${\rm;}
in particular, $f\in M \Longrightarrow g\in M$ and $f\in M_c \Longrightarrow g\in M_c$.
Dually, if the square is a pushout then $\coker g=(\coker f)\beta$ and $g=\coker\zeta$
implies $f=\coker(\alpha\zeta)${\rm;} in particular, $g\in P \Longrightarrow f\in P$ and
$g\in P_c \Longrightarrow f\in P_c$.
\end{lemma}

An additive category meeting Axiom~\ref{1} is abelian if and only
if $\overline\alpha$ is an isomorphism for every
$\alpha$. Consider the following axiom.

\begin{axiom}\label{a2}
For every morphism
$\alpha$,
$\overline\alpha$
is a bimorphism, i.e., a monomorphism and an epimorphism.
\end{axiom}

We write $\alpha\|\beta$ if the sequence
$\cdot\stackrel{\alpha}{\longrightarrow}
\cdot\stackrel{\beta}{\longrightarrow}\cdot$
is exact, that is, $\im\alpha=\ker\beta$ (which, in a category meeting
Axioms~1 and~2, is equivalent to $\coker\alpha=\coim\beta$).

\begin{lemma}[\protect\cite{KK1}]
The following assertions hold in an additive category satisfying
Axioms~{\rm\ref{a1}} and~{\rm\ref{a2}:}

{\rm(1)} if $gf\in M_c$ then $f\in M_c${\rm;} if $gf\in P_c$ then
$g\in P_c${\rm;}

{\rm(2)} if $f,g\in M_c$ and $fg$ is defined then $fg\in M_c$, if
$f,g\in P_c$ and $fg$ is defined then $fg\in P_c${\rm;}

{\rm(3)} if $fg\in O_c$ and $f\in M$ then $g\in O_c$,
if $fg\in O_c$ and $g\in P$ then $f\in O_c$.
\end{lemma}

It is well known (see, for example,~\protect\cite{PP}), that every abelian
category satisfies the following two axioms dual to one another.

\begin{axiom}\label{a3}
If~{\rm(1)} is a pullback then $f\in P_c\Longrightarrow g\in P_c$.
\end{axiom}

\begin{axiom}\label{a4}
If~{\rm(1)} is a pushout then $g\in M_c\Longrightarrow f\in M_c$.
\end{axiom}

An additive category satisfying Axioms~\ref{a1}, \ref{a3}, and~\ref{a4}, is
called {\it quasi-abelian}. Such categories are also known
as (\textit{Ra{\u\i}kov})-\textit{semiabelian} (the original name, proposed 
by Ra{\u\i}kov in~\protect\cite{Ra} and used in the Russian tradition; now, however,
the term {\it semi-abelian category} is involved in a quite
different context~\protect\cite{JMT}) or
\textit{almost abelian}~\protect\cite{R2}. As follows from
Theorem~1 of~\protect\cite{KC}, each quasi-abelian category meets
Axiom~\ref{a2}.

Given an arbitrary commutative square~(\ref{1}), denote by
$\widehat g:\Ker\alpha\to\Ker\beta$
the morphism defined by the equality
$g(\ker\alpha)=(\ker\beta)\widehat g$
and by
$\widehat f:\Coker\alpha\to\Coker\beta$
the morphism defined by the condition
$\widehat f(\coker\alpha)=(\coker\beta)f$.

From now on, unless otherwise specified, the ambient category
$\mathcal{A}$  is assumed quasi-abelian.

Lemmas~5 and~6 of~\protect\cite{KK1} yield the following assertion.

\begin{lemma}[\protect\cite{KK1}]\label{l3}
Suppose that square~{\rm(\ref{1})} is a pullback. If $\beta\in O_c$ then
$\alpha\in O_c$ and $\widehat f\in M$.

Dually, if~{\rm(\ref{1})} is a pushout and $\alpha\in O_c$ then $\beta\in O_c$ and
$\widehat g\in P$.
\end{lemma}

\begin{lemma}[Composition Lemma]\label{l4}
Suppose that the composition $gf$ of two morphisms $f$ and $g$ is defined.
Then there exists a semiexact sequence
\begin{multline}\label{4}
0\to\Ker f\stackrel{\varphi}{\to} \Ker(gf) \stackrel{\psi}{\to}
\Ker g \stackrel{\chi}{\to} \Coker f \\
\stackrel{\lambda}{\to}
\Coker(gf) \stackrel{\omega}{\to} \Coker g \to 0
\end{multline}
which is exact at~$\Ker f$, $\Ker(gf)$, $\Coker(gf)$, and $\Coker g$;
moreover, $\varphi$ and $\omega$ are strict. Furthermore, if $f\in O_c$
then~{\rm(\ref{4})} is exact at~$\Ker g$ and $\psi\in O_c$; if $g\in O_c$
then~{\rm(\ref{4})} is exact at~$\Coker f$ and $\chi\in O_c$.
\end{lemma}

\begin{proof}
As in an~abelian category, we define $\varphi$, $\psi$,
$\chi$, $\lambda$, and $\omega$ by the equalities $\ker f=(\ker(gf))\varphi$,
$f(\ker(gf))=(\ker g)\psi$, $\chi=(\coker f)(\ker g)$,
$(\coker(gf))g=\lambda(\coker g)$, and $\coker g=\omega(\coker(gf))$.
Then it is standard (and easy) that sequence~(\ref{4}) thus obtained is
semiexact, $\varphi=\ker\psi$, and $\omega=\coker\lambda$. Furthermore,
it is easy to check that the square
\begin{equation}\label{5}
\bfig
\morphism(0,200)<600,0>[\Ker(gf)`\cdot;\ker(gf)]
\morphism(0,200)<0,-400>[\Ker(gf)`\Ker g;\psi]
\morphism(0,-200)<600,0>[\Ker g`\cdot;\ker g]
\morphism(600,200)<0,-400>[\cdot`\cdot;f]
\efig
\end{equation}
is a~pullback.

Suppose that $f$ is strict. Applying Lemma~3 to pullback~(\ref{5}), we obtain
that $\psi\in O_c$ and the morphism $l$ defined by the equality
$l(\coker\psi)=(\coker f)(\ker g)\,\,(=\chi)$ is monic. Thus,
$\im\psi=\ker\chi$, which proves the exactness at~$\Ker g$. By duality,
we infer that $\lambda\in O_c$ and~(\ref{4}) is exact at~$\Coker f$.
The lemma is proved.
\end{proof}

\section{Lambek Invariants}\label{sec2}

Given a~commutative square~(\ref{1}), consider the pullback
\begin{equation}\label{6}
\bfig
\morphism(0,200)<400,0>[I`\Img f;k]
\morphism(0,200)<0,-400>[I`\Img\beta;l]
\morphism(0,-200)<400,0>[\Img\beta`B;\im\beta]
\morphism(400,200)<0,-400>[\Img f`B;\im f]
\efig
\end{equation}
Easily, there are morphisms $k':\Img(f\alpha)\to\Img f$ and
$l':\Img(f\alpha)\to\Img\beta$ with
$\im(f\alpha)=(\im f)k'=(\im\beta)l'$. Since~(\ref{6}) is a~pullback,
it follows that there is a~unique morphism $\rho:\Img(f\alpha)\to I$ such that
$k'=k\rho$ and $l'=l\rho$. We put $\Img S=\Coker\rho$. If we denote by
$\Phi$ the epimorphism $\widetilde{f\alpha}$ then, obviously,
$\Img S=\Coker(\rho\Phi)$.

Now, let $\mu:\Ker g\to\Ker(f\alpha)$ and $\nu:\Ker\alpha\to\Ker(f\alpha)$
be the natural inclusions. They form a~morphism
$\langle\mu,\nu\rangle:\Ker g\oplus\Ker\alpha\to\Ker(f\alpha)$. We put
$\Ker S=\Coker\langle\mu,\nu\rangle$. Alternatively, $\Ker S$ can be described
as follows (see, for example, \protect\cite{No1}). Consider the pushout
\begin{equation*}
\bfig
\morphism(0,200)<600,0>[C`\Coim\alpha;\coim\alpha]
\morphism(0,200)<0,-400>[C`\Coim g;\coim g]
\morphism(0,-200)<600,0>[\Coim g`J;i]
\morphism(600,200)<0,-400>[\Coim\alpha`J; j]
\efig
\end{equation*}
There is a unique morphism $\sigma:L\to B$ such that
$\sigma j=f(\im\alpha)\bar\alpha$ and $\sigma i=\beta(\im g)\bar g$.
Then $\Ker S$ is naturally isomorphic with $\Ker\sigma$. Thus, $\Img S$
and $\Ker S$ are dual notions.

In what follows, we endow all the morphisms and objects introduced above for
a~commutative square~$S$ with the subscript $S$ when it becomes necessary
to distinguish the corresponding morphisms of different squares.

The condition $\Img S=0$ ($\Ker S=0$) is fulfilled for an important class
of pullbacks in a~quasi-abelian category. Namely, the following
assertion holds.

\begin{them}\label{t1}
Suppose that square~{\rm(\ref{1})} is a~pullback with $\beta$ and $f$
strict. Then $\Img S=0$. If~{\rm(\ref{1})} is a~pushout with $\alpha$
and $g$ strict then $\Ker S=0$.
\end{them}

\begin{proof}
Consider the commutative diagram
$$
\bfig
\morphism(0,500)<500,0>[F`\cdot;v_1]
\morphism(500,500)<500,0>[\cdot`D;w_1]
\morphism(0,500)<0,-500>[F`\cdot;v_2]
\morphism(500,500)<0,-500>[\cdot`I;v_0]
\morphism(1000,500)<0,-500>[D`\Img f;\coim f]
\morphism(0,0)<500,0>[\cdot`I;]
\morphism(500,0)<500,0>[I`\Img f;k]
\morphism(0,0)<0,-500>[\cdot`A;w_2]
\morphism(500,0)<0,-500>[I`\Img\beta;l]
\morphism(1000,0)<0,-500>[\Img f`B,;\im f]
\morphism(0,-500)<500,0>[A`\Img\beta;\coim\beta]
\morphism(500,-500)<500,0>[\Img\beta`B,;\im\beta]
\efig
$$
where all the four squares are pullbacks. Then the ``resulting'' square
is a~pullback, too (see, for example, \protect\cite{EH1}, Proposition~2.10).
Thus, up to an~isomorphism, we have $C=F$, $w_1 v_1=\alpha$, and
$w_2 v_2=g$. Since $w_1,w_2\in~M_c$ and $v_1,v_2\in P_c$, by Lemma~\ref{l1}(3)
it follows that $w_1=\im\alpha$, $v_1=\coim\alpha$, $w_2=\im g$, and
$v_2=\coim g$. Therefore, $\im(f\alpha)=\im((\im f)k v_0 v_1)=(\im f)k$, and
hence $I=\Img(f\alpha)$, which implies $\Img S=~0$.

The second assertion is proved by duality.

The theorem is proved.
\end{proof}

\begin{remk}
{\rm
By Lemma~\ref{l3}, if square~(\ref{1}) is a~pullback with $\beta\in O_c$
($f\in O_c$) then $\alpha\in O_c$ ($g\in O_c$). This means that Theorem~\ref{t1}
applies to ``strict'' pullbacks. However, it fails to hold for ``nonstrict''
pullbacks, which is demonstrated by the following example. Consider
the category  $\mathcal{B}an$ of Banach spaces and bounded linear operators.
Let $A$ and $B$ be infinite-dimensional Banach spaces and let $\beta:A\to B$
be a linear operator with dense range $R(\beta)\ne B$ (and so
$\beta\not\in O_c$!). Put $D=\mathbb{R}$
and suppose that $f:D\to B$ is injective and $R(f)\cap R(\beta)=0$.
Form a~pullback $f\alpha=\beta g$. For a~morphism $L:X\to Y$ in
$\mathcal{B}an$, $\Img L$ is the closure $\overline{R(L)}$ of its range $R(L)$. It is easy
to see that $\alpha=0$ and hence $\Img(f\alpha)=0$. However, in this case,
$I=\overline{R(\beta)}\cap\overline{R(f)}\cong\mathbb{R}\ne 0$. Thus,
$\Img S\cong\mathbb{R}$. }
\end{remk}

\begin{remk}
{\rm
The set of commutative squares $S$ with $\Img S=0$
($\Ker S=0$) is not reduced to ``strict'' pullbacks (pushouts). As observed by
Hilton (see~\protect\cite{Hi}, Proposition~2.4) and is easily checked, each
composition $h=gf$ yields two commutative squares $\Delta':\, h(\id)=gf$
and $\Delta'':\, (\id)h=gf$ such that $\Img\Delta'=0$ and $\Ker\Delta''=0$.
Obviously, $\Delta'$ is a~pullback if and only if $g$ is monic
(similarly, $\Delta''$ is a~pushout if and only if $f$ is epic). Hence,
a~commutative square $S$ need not be a~pullback (pushout) to have $\Img S=0$
($\Ker S=0$).  }
\end{remk}

As noted in the introduction, for a~sequence of the form of~(\ref{2}) with
exact rows, $\Ker S$ and $\Img T$ are known to be naturally isomorphic
(see \protect\cite{Le} or \protect\cite{No1}) in an exact category.
For this to hold in a~quasi-abelian category, one must
have $\Img b=\Coim b$, that is, $b$ must be strict. On the same assumption,
we can use Nomura's construction of $\Lambda:\Img S\to\Ker T$ for a~diagram
of the kind of~(\ref{2}) with semiexact rows. Recall that $\Lambda$ is
characterized by the equality
$(\ker\sigma_T)\Lambda(\coker\rho_S)=i_T k_S$~\protect\cite{No1}.

When the rows in~(2) are exact, we can still construct a~canonical morphism
$\zeta:\Ker T\to \Img S$. Of course, $\zeta=\Lambda^{-1}$ if $\Lambda$ exists.
Namely, we have

\begin{them}\label{t2}
Suppose that in~{\rm(\ref{2})} the rows are exact. Then there exist
unique morphisms $\xi:\Ker(g'b)\to I_S$ and $\zeta:\Ker T\to\Img S$ such
that
$$
(\coker\rho_s)\xi=\zeta\coker\langle\mu_T,\nu_T\rangle.
$$
\end{them}

\begin{proof}
Obviously, $g'b(\ker(g'b))=0$, which implies that there exists
a~unique morphism $y$ with $b\ker(g'b))=(\ker g')y_0=(\im f')y_0$.
Since~(\ref{6}) is a~pullback, it follows that there exists a~unique morphism
$\xi:\Ker(g'b)\to\Img S$ such that $\tilde b(\ker(g'b))=k_S\xi$ and
$y=l_S\xi$. We have
$$
k_S\xi\mu_T=\tilde b(\ker(g'b))\mu_T=\tilde b(\ker b)=0,
$$
whence $\xi\mu_T=0$ because $k_S$ is monic. Now, denote by $\gamma=\gamma_S$
the unique morphism for which $\im(bf)\gamma=b(\im f)$ ($=b(\ker g)$ by
the exactness of the upper row in~(\ref{2})). We infer
\begin{multline*}
(\im b)k_S\rho_S \gamma_S\tilde f=(\im(bf))\gamma_S\tilde f
=b(\im f)\tilde f=bf \\
=(\im b)\tilde b(\ker(g'b))\nu_T\tilde f
=(\im b)k_S\xi\nu_T\tilde f.
\end{multline*}
Since $(\im\varphi)k\in M$ and $\tilde f\in P$, it follows that
$\xi\nu_T=\rho_S\gamma_S$. Hence
$(\coker\rho_S)\xi\nu_T=(\coker\rho_S)\xi\langle\mu_T,\nu_T\rangle=0$.
Therefore, there exists a~unique morphism
$\zeta:\Coker\langle\mu,\nu\rangle\to\Coker\rho_S$ such that
$$
(\coker\rho_S)\xi=\zeta\coker\langle\mu_T,\nu_T\rangle.
$$
The theorem is proved.
\end{proof}

As a corollary to Theorem~\ref{t2}, we obtain Lambek's isomorphism, established for
exact categories in~\protect\cite{La,Le,No1}, which, in our case, holds
under the extra assumption that $b\in O_c$. Note that, in view of the exactness
properties of the $\Ker$-$\Coker$-sequence in a quasi-abelian
category proved in~\protect\cite{KK1}, Nomura's proof of Lambek's isomorphism
in~\protect\cite{No1} is carried over to our situation literally. However,
here we prefer to show how $\zeta$ becomes an isomorphism if $b$ is strict.

\begin{corl}
If, under the conditions of Theorem~{\rm\ref{t2}},
$b\in O_c$ then $\zeta$ is an isomorphism.
\end{corl}

\begin{proof}
First, observe that the square
\begin{equation}\label{7}
\bfig
\morphism(0,200)<500,0>[\Ker(g'b)`B;\ker(g'b)]
\morphism(0,200)<0,-400>[\Ker(g'b)`I;\xi]
\morphism(0,-200)<500,0>[I`\Img b;k_S]
\morphism(500,200)<0,-400>[B`\Img b; \tilde b]
\efig
\end{equation}
is a pullback.

Indeed, suppose that morphisms $x_1$ and $x_2$ are such that
$k_S x_1=\tilde b x_2$. Then
$$
g'b x_2=g'(\im b)\tilde b x_2 =g'(\im b)k_S\xi \\
= g'(\im f')l_S\xi =g'(\ker g')l_S\xi=0.
$$
Therefore, there exists a unique morphism $x$ with $x_2=(\ker(g'b))x$.
We now prove that $x_1=\xi x$. We have
\begin{multline*}
(\im b)k_S\xi x=(\im f')l_S\xi x = (\ker g')l_S\xi x = b(\ker(g'b))x \\
=(\im b)\tilde b(\ker(g'b))x=(\im b)\tilde b x_2=(\im b)k_S x_1,
\end{multline*}
whence, by the fact that $(\im b)k_S$ is monic, we see that $\xi x= x_1$.
Thus, we have demonstrated that (\ref{7}) is a~pullback.

Since $\ker b=\ker\tilde b=(\ker(g'b))\mu_T$, $b\in O_c$, and (7) is a pullback, from
Lemma~\ref{l1}(4) and Axiom~3 it follows that $\mu_T\|\xi$. Obviously, we have
$(\coker\langle\mu_T,\nu_T\rangle) \mu_T=~0$,
and so there exists a unique morphism
$\tau:I\to\Coker\langle\mu_T,\nu_T\rangle$ such that
$\coker\langle\mu_T,\nu_T\rangle=\tau\xi$.
We have $\zeta\tau\xi=(\coker\rho_S)\xi$, and the relation $\xi\in P_c$
yields $\zeta\tau=\coker\rho_S$. Furthermore,
\begin{equation}\label{8}
\tau\rho_S \gamma_S=\tau\xi\nu_T=(\coker\langle\mu_T,\nu_T\rangle)\nu_T=0.
\end{equation}
Since
$\gamma_s\tilde f=\Phi$, it follows that $\gamma_S$ is epic and so (\ref{8})
implies that $\tau\rho_S=0$. Thus there is a unique morphism
$\Lambda_0:\Coker\rho\to\Coker\langle\mu_T,\nu_T\rangle$ with the property
$\tau=\Lambda_0(\coker\rho)$. Easily, $\zeta\Lambda_0$ and $\Lambda_0\zeta$
are identities and, therefore, $\zeta$ and $\Lambda_0$ are mutually inverse
isomorphisms. This finishes the proof of the corollary.
\end{proof}

It can be proved that, up to the identification
$\Ker T\cong \Ker\sigma_T$, $\Lambda_0$ is Nomura's morphism $\Lambda$.

We now pass to the more general case of a commutative diagram of the form
of~(\ref{2}) with semiexact rows.

In the case of an~exact category, Nomura constructed exact sequence~(\ref{3}).
However, an analysis of the proof of the exactness of~(\ref{3})
in~\protect\cite{No1} (based on the Composition Lemma, cf. Lemma~\ref{4})
shows that, in the quasi-abelian case, many morphisms must
be assumed strict so that all morphisms in~(\ref{3}) can be defined.
We prove the following quasi-abelian version of Corollary~$A_2$
of \protect\cite{No1}.

\begin{them}
Suppose that in diagram~{\rm(\ref{2})} the rows are semiexact. The following
asserions hold.

{\rm (1)} If the sequence $A'\rightarrow B'\rightarrow C'$ is exact
and $b\in~M_c$ then there exists a canonical morphism
$\theta:H(A\rightarrow B\rightarrow C)\rightarrow \Img S$ such that
the sequence
$$
0 \to H(A\rightarrow B\rightarrow C) \stackrel\theta\to\Img S
\stackrel\Lambda\to \Ker T\to 0
$$
is exact.

{\rm (2)} If the sequence $A\rightarrow B\rightarrow C$ is exact and
$b\in P_c$ then there exists a~canonical morphism
$\varkappa:\Ker T\rightarrow H(A'\rightarrow B'\rightarrow C')$
such that the sequence
$$
0 \to \Img S \stackrel\Lambda\to \Ker T \stackrel\varkappa\to
H(A'\rightarrow B'\rightarrow C')\to 0
$$
is exact.
\end{them}

\begin{proof}
We prove only item~1 because item~2 is obtained from it by duality.

By definition, the homology object
$H(A\rightarrow B\rightarrow C)$ is the cokernel of a~unique
morphism $\varepsilon$ such that $\im f=(\ker g)\varepsilon$.
Consequently, $(\coker\rho_S)\xi\nu_T\varepsilon=0$ and, therefore,
there exists a unique morphism $\theta$ with
$(\coker\rho_S)\xi\nu_T=\theta(\coker\varepsilon)$.
Repeating the argument of the proof of Theorem~\ref{2} almost literally,
we see that $\rho_S \gamma_S= \xi\nu_T\varepsilon$. Furthermore, since
$b(\im f)=(\im(bf))\gamma_S$ and $b$ is a kernel, it follows that $\gamma_S$
is an isomorphism.  In addition, $\xi$ is an isomorphism, too. Indeed,
as above, $\xi$ is a~part of pullback~(\ref{7}), which implies that $\xi\in P_c$
and $(\ker(g'b))(\ker\xi)=\ker\tilde b=(\ker(g'b))\mu_T=0$. Thus $\mu_T=0$
and hence $\xi$ is in fact an isomorphism. Thus we may write
$\rho_S=\nu_T\varepsilon$. Since we thus obtain a pullback
$\rho_S\id = \nu_T\varepsilon$ and $\rho\in O_c$, the morphism of the
cokernels $\theta:\Coker\varepsilon\to\Coker\rho_S$ is monic. Thus
we see the exactness at $H(A\rightarrow B\rightarrow C)$.

Furthermore, since
$$
\Lambda\theta(\coker\varepsilon)=\Lambda(\coker\rho_S)\xi\nu_T
=(\coker\langle\mu_t,\nu_T\rangle)\nu_T=~0,
$$
we infer $\Lambda\theta=0$. Now, take a~morphism $y$
with $y\theta=0$. Then $y(\coker\rho_S)\nu=y\theta(\coker\varepsilon)=0$
and, obviously, $y(\coker\rho_S)\mu_T=0$. Hence, there exists a~unique
morphism $v$ with $y(\coker\rho_S)=v(\coker\langle\mu_T,\nu_T\rangle)
=v\Lambda(\coker\rho_S)$. Since $\coker\rho_S$  is epic, it follows that
$y=v\Lambda$. Thus, $\Lambda=\coker\theta$ and so we have the exactness
at $\Img S$.

Theorem~\ref{3} is proved.
\end{proof}

We now prove another assertion about a~diagram of commutative squares
(cf. Proposition~2.7 in~\protect\cite{Hi}).

\begin{them}
Suppose that, in a~commutative diagram
\begin{equation}\label{9}
\bfig
\morphism(0,400)<400,0>[A_1`B_1;\theta_1]
\morphism(400,400)<400,0>[B_1`C_1;\theta_2]
\morphism(0,400)<0,-400>[A_1`A_2;\alpha_1]
\morphism(400,400)<0,-400>[B_1`B_2;\beta_1]
\morphism(800,400)<0,-400>[C_1`C_2;\gamma_1]
\morphism(0,0)<400,0>[A_2`B_2;\varphi_1]
\morphism(400,0)<400,0>[B_2`C_2;\varphi_2]
\place(200,200)[S]
\place(600,200)[T]
\morphism(0,0)<0,-400>[A_2`A_3;\alpha_2]
\morphism(400,0)<0,-400>[B_2`B_3;\beta_2]
\morphism(800,0)<0,-400>[C_2`C_3;\gamma_2]
\morphism(0,-400)<400,0>[A_3`B_3;\psi_1]
\morphism(400,-400)<400,0>[B_3`C_3;\psi_2]
\place(200,-200)[U]
\place(600,-200)[V]
\efig,
\end{equation}
the first column is exact at~$A_2$, the third, at~$C_2$, and
the second row is exact at~$B_2$, $\Img T=0$, $\Ker U=0$,
$\beta_2\beta_1=0$, $\varphi_1\in O_c$, and $\varphi_2\beta_1\in O_c$.
Then the second column is exact at~$B_2$.
\end{them}

\begin{proof}
Take a~morphism $x:X\to B_2$ such that $\beta_2 x=0$. We may assume
that $x=\im x$. We have $\gamma_2\varphi_2 x=0$; therefore, there exists
a~unique morphism $y$ such that $\varphi_2 x=(\im\gamma_1)y$. Since
$(\im\varphi_2)\tilde\varphi_2 x=(\im\gamma_1)y$ and $\Img T=0$, it follows
that there is a~unique morphism $\xi:X\to \Img(\varphi_2\beta_1)$ with
the properties $\tilde\varphi_2 x=l_T\xi$ and $y=k_T\xi$. Thus,
$\varphi_2 x=(\im\varphi_2) l\xi=\im(\varphi_2\beta_1)\xi$. Define $\omega$
by the equality $\varphi_2\beta_1=\im(\varphi_2\beta_1)\omega$. Then
$\omega\in P_c$. Consider a~pullback $\xi\omega_0=\omega\xi_0$.
We have $\im(\varphi_2\beta_1)\xi\omega_0=\im(\varphi_2\beta_1)\omega\xi_0
=\varphi_2\beta_1\xi_0$. Thus, $\varphi_2(x\omega_0-\beta_1\xi_0)=0$,
whence we deduce the existence of a~unique morphism $\xi_1$ such that
$x\omega_0-\beta_1\xi_0=(\ker\varphi_2)\xi_1=(\im\varphi_1)\xi_1$.
Let $\xi_0 p_0=\tilde\varphi_1\xi_1$ be a~pullback. Then
$$
0=\beta_2 x\omega_0=\beta_2(\im\varphi_1)\xi_0 p_0
=\beta_2(\im\varphi_1)\tilde\varphi_1\xi_1=\beta_2\varphi_1\xi_1.
$$
Since $\Ker U=0$, it follows that
$\Ker(\beta_2\varphi_1)\cong \Ker\varphi_1\oplus\Ker\alpha_2$.
Consequently, $\xi_1=(\ker\varphi_1)t_1+(\ker\alpha_2)t_2
=(\ker\varphi_1)t_1+(\im\alpha_1)t_2$ for some $t_1$ and $t_2$.
Furthermore, there exists a~unique morphism $u$ with
$\varphi_1(\im\alpha_1)=(\im\beta_1)u$. We infer
$$
x\omega_0 p_0=\beta_1\xi_0 p_0+\varphi_1(\im\alpha_1)t_2
=\beta_1\xi_0 p_0 + (\im\beta_1)u t_2
=(\im\beta_1)(\tilde\beta_1\xi_0 p_0+ut_2).
$$
Thus, $x\omega_0 p_0=(\im\beta_1)v$,
i.~e.,
$x\omega_0 p_0=(\im\beta_1)(\im v)\bar v(\coim v)$.
The hypothesis implies that $\omega_0 p_0\in P_c$. Therefore,
$x=(\im\beta_1)(\im v)$, which means that $\im\beta_1=\ker\beta_2$.

Theorem~\ref{4} is proved.
\end{proof}

For abelian categories, Theorem~\ref{4} was proved by Hilton
(see~\protect\cite{Hi}, Proposition~2.7) and served as a~key ingredient
in the proof of the~main theorem in~\protect\cite{Hi} on the~exactness of
a~system of interlocking exact sequences. In the quasi-abelian
case, we have to add some strictness conditions to Hilton's Proposition~2.7.
Unfortunately, applying Theorem~\ref{4} to interlocking sequences
(and thus to spectral sequences) is possible only if we assume all
the morphisms strict. We dealt with spectral sequences by considering
exact couples in quasi-abelian categories in a~separate
paper~\protect\cite{K}.

\textbf{Acknowledgments.}
The author acknowledges the financial support of
a postdoctoral (research) fellowship from NATO,
INTAS Grant 03--51--3251, and the State Maintenance Program for Leading
Scientific Schools of the Russian Federation (Grant~ NSh~311.2003.1).

This paper was begun in 2003, when the author was
a NATO postdoctoral fellow at the Universit\'e de Lille~I. It is a pleasure
for him to thank this university and especially Prof. Leonid Potyagailo for
the hospitality he enjoyed during his nine-month stay in Lille.

All diagrams in the paper were made using Prof. Michael Barr's
\texttt{diagxy} macro package.


\begin{thebibliography}{25}

\bibitem{BD}
Bucur~I and Deleanu~A.:
\newblock{\em Introduction to the Theory of Categories and Functors,}
\newblock{Pure and Applied Mathematics, XIX, Interscience Publication John
Wiley \& Sons, Ltd., London-New York-Sydney, 1968.}

\bibitem{EH1}
Eckmann~B. and Hilton~P.~J.: Exact couples in an abelian category,
{\em J. Algebra} \textbf{3} (1966), 38--87.

\bibitem{GK1}
Glotko~N.~V. and Kuz$'$minov~V.~I.:
On the cohomology sequence in a semiabelian category (Russian),
{\em Sib. Mat. Zh.} \textbf{43} (2002), no. 1, 41--50; English translation in:
{\em Sib. Math. J.} \textbf{43} (2002), no. 1, 28--35.

\bibitem{GK2}
Glotko~N.~V. and Kuz$'$minov~V.~I.:
On reflective subcategories of quasiabelian categories,
{\em Sib. \`Elektron. Mat. Izv.} \textbf{2} (2005), 68-78.
(Russian); http://semr.math.nsc.ru/V2/v2p68-78.pdf.

\bibitem{Hi}
Hilton~P.~J.: On systems of interlocking exact sequences,
{\em Fundam. Math.} \textbf{61} (1967), 111-119.

\bibitem{JMT}
Janelidze~G., M\'arki~L., and Tholen~W.:  Semi-abelian categories,
Category theory 1999 (Coimbra), {\em J. Pure Appl. Algebra} \textbf{168}
(2002), no. 2-3, 367--386.

\bibitem{K}
Kopylov Ya.~A.: Exact couples in a~Raikov-semiabelian category,
{\em Cah. Topol. G\'eom. Diff\'er. Cat\'eg.} \textbf{45} (2004), no.~3, 162--178.

\bibitem{KK1}
Kopylov~Ya.~A. and Kuz$'$minov~V.~I.: On the Ker-Coker sequence in a
semiabelian category, {\em Sib. Mat. Zh.} \textbf{41} (2000), no. 3, 615--624;
English translation in: {\em Sib. Math. J.} \textbf{41} (2000), no. 3, 509--517.

\bibitem{KK2}
Kopylov~Ya.~A. and Kuz$'$minov V.~I.: Exactness of the cohomology sequence
corresponding to a short exact sequence of complexes in a semiabelian
category, {\em Sib. Adv. Math.} \textbf{13} (2003), no.~3, 72--80.

\bibitem{KC}
Kuz$'$minov~V.~I. and Cherevikin~A.~Yu.: Semiabelian categories,
{\em Sib. Mat. Zh.} {\textbf 13} (1972), no.~6, 1284--1294; English translation in:
{\em Sib. Math. J.} \textbf{13} (1972), no.~6, 895--902.

\bibitem{La}
Lambek J.:
Goursat's theorem and homological algebra, {\em Can. Math. Bull.}
\textbf{7}, (1964) 597-608.

\bibitem{Le}
Leicht~J.~B.:
Axiomatic proof of J. Lambek's homological theorem,
{\em Can. Math. Bull.} \textbf{7} (1964) 609-613.

\bibitem{No1}
Nomura~Y.:
An exact sequence generalizing a theorem of Lambek,
{\em Arch. Math.} \textbf{22} (1971), 467-478.

\bibitem{No2}
Nomura~Y.:
Induced morphisms for Lambek invariants of commutative squares,
{\em Manuscr. Math.} \textbf{4} (1971), 263-275.


\bibitem{PP}
Popescu~N. and Popescu~L.:
\newblock{\em Theory of categories},
{\newblock Bucuresti: Editura Academiei.
Alphen aan den Rijn: Sijthoff \& Noordhoff International Publishers, 1979.}

\bibitem{P1}
Prosmans~F.: Derived projective limits of topological abelian groups,
{\em J. Funct. Anal.} \textbf{162} (1999), no. 1, 135--177.

\bibitem{P2}
Prosmans~F.: Derived limits in quasi-abelian categories,
{\em Bull. Soc. Roy. Sci. Li\`ege} \textbf{68} (1999),
no. 5-6, 335--401.

\bibitem{P3}
Prosmans~F.: Derived categories for functional analysis,
{\em Publ. Res. Inst. Math. Sci.} \textbf{36} (2000), no. 1, 19--83.

\bibitem{Ra}
Ra\u\i kov~D.~A., Semiabelian categories, {\em Dokl. Akad. Nauk SSSR}
{\textbf 188} (1969), 1006--1009; English translation in
{\em Soviet Math. Dokl.} \textbf{10} (1969), 1242-1245.

\bibitem{R1}
Rump~W.: $*$-modules, tilting, and almost abelian categories,
{\em Comm. Algebra} \textbf{29} (2001), no. 8,
3293--3325; Erratum (Misprints generated via electronic editing):
{\em Comm. Algebra} \textbf{30} (2002), 3567--3568.

\bibitem{R2}
Rump~W.: Almost abelian categories,
{\em Cah. Topol. G\'eom. Diff\'er. Cat\'eg.}
\textbf{42} (2001), no. 3, 163--225.

\bibitem{R3}
Rump~W.:
Categories of lattices, and their global structure in terms of
almost split sequences, {\em Algebra Discrete Math.} (2004), no. 1,
87--111.

\bibitem{S}
Schneiders~J.-P.: {\em Quasi-Abelian Categories and Sheaves}, M\'em. Soc. Math. Fr.
(N.S.) (1999), no. 76.

\bibitem{Sc}
Succi Cruciani~R.: Sulle categorie quasi abeliane,
{\em Rev. Roumaine Math. Pures Appl.} \textbf{18} (1973), 105--119.

\end{thebibliography}
\end{document}